\documentclass[onecolumn,12pt]{elsarticle}

\usepackage{amssymb}

\journal{arXiv.org}

\newtheorem{theorem}{Theorem}[section]

\newtheorem{remark}[theorem]{Remark}
\newtheorem{example}[theorem]{Example}

\newfont{\cirilrm}{wncyr10 scaled 1100}
\newfont{\cirilit}{wncyi10 scaled 1100}
\newfont{\cirilsc}{wncysc10 scaled 1100}

\newcommand{\ja}{\symbol{'37}}

\begin{document}

\begin{frontmatter}

\title{
Differential transcendency in {the} theory of linear differential systems with constant coefficients
}

\medskip

\author{
Branko Male\v sevi\' c\mbox{\small ${\,}^{a}$},
Dragana Todori\' c\mbox{\small ${\,}^{b}$},
Ivana Jovovi\' c\mbox{\small ${\,}^{a}$},
Sonja Telebakovi\' c\mbox{\small ${\,}^{b}$}
}

\medskip

\address[etf]{Faculty of Electrical Engineering, University of Belgrade, \\
Bulevar kralja Aleksandra 73, 11000 Belgrade, Serbia}
\address[matf]{Faculty of Mathematics, University of  Belgrade, \\
Studentski trg 16, 11000 Belgrade, Serbia}

\medskip

\begin{abstract}
In this paper we consider a reduction of a non-homogeneous linear system of first order operator equations to a~totally reduced system.
Obtained results are applied to Cauchy problem for linear differential systems with constant coefficients and to the question
of differential transcendency.
\end{abstract}

\medskip

\end{frontmatter}

\renewcommand{\thefootnote}{}

\footnotetext{
{\sl Email addresses}$\,:$

\quad
Branko Male\v sevi\' c $<${\sl malesevic@etf.rs}$>\,$,
Dragana Todori\' c $<${\sl draganat@math.rs}$>\,$,

\quad
Ivana Jovovi\' c $<${\sl ivana@etf.rs}$>\,$,
Sonja Telebakovi\' c $<${\sl sonjat@math.rs}$>$}

\bigskip

\section{Introduction}

\bigskip

Linear systems with constant coefficients are considered in various fields
{\big (}$\!\;$see \cite{Fortmann--Hitz--77, Martínez-Guerra--González-Galan--Luviano-Juárez--Cruz-Victoria--07},
\cite{Cruz-Victoria--Gonzalez-Sanchez--08}-\cite{MalesevicTodoricJovovicTelebakovic2010} $\!${\big )}.
In our paper \cite{MalesevicTodoricJovovicTelebakovic2010} we use the rational canonical form
and a certain sum of principal minors to reduce a linear system of first order ope\-rator equations
with constant coefficients to an equivalent, so called partially reduced, system. In this paper
we obtain more general results regarding sums of principal minors and a new type of reduction.
The obtained formulae of reduction allow some new considerations in connection with Cauchy problem
for linear differential systems with constant coefficients and {in connection with} the differential
transcendency of the solution coordinates.

\medskip

\section{Notation}

\bigskip

Let us recall some notation. Let $K$ be a field and let $B \!\in\! K^{n \times n}$
be an $n$-square matrix. We denote by$:$
\begin{description}
\item[($\mbox{\large $\delta$}_{k}$)]
$\mbox{\large $\delta$}_{k} \!=\!\mbox{\large $\delta$}_{k}(B)$
the sum of its principal minors of order $k$ ($1 \leq k \leq n$),
\item[($\mbox{\large $\delta$}^{i}_{k}$)]
$\mbox{\large $\delta$}^{i}_{k}\!=\!\mbox{\large
$\delta$}^{i}_{k}(B)$ the sum of its principal minors of order $k$
containing~$i^{th}$ column $(1 \leq i, k \leq n)$.
\end{description}

Let $v_{1}, \ldots, v_{n}$ be elements of $K$. We write $B^{i}(v_{1}, \ldots, v_{n})\!\in\! K^{n \times n}$ for the matrix
obtained by substituting $\vec{v}\!=\![ v_{1}\;\ldots\; v_{n}]^{T}$ in place of $i^{th}$ column of $B$. Furthermore, it is
convenient to use
\begin{equation}
\label{lepotani}
\mbox{\large $\delta$}^{i}_{k}(B; \vec{v})
=
\mbox{\large $\delta$}^{i}_{k}(B; v_{1}, \ldots, v_{n})
=
\mbox{\large $\delta$}^{i}_{k}(B^i (v_{1}, \ldots, v_{n})),
\end{equation}
and the corresponding vector
\begin{equation}
\mbox{\large $\vec{\delta}$}_{k}{\big (} B; \vec{v}
{\big )} \!=\! [\;\mbox{\large $\delta$}_{k}^{1}{\big (} B;
\vec{v} {\big )} \;\;\ldots\;\;\mbox{\large $\delta$}_{k}^{n}{\big
(} B; \vec{v} {\big )}\;]^T.
\end{equation}
The characteristic polynomial $\Delta_{B}(\lambda)$ of the matrix $B \!\in\!K^{n \times n}$ has the form
$$
\Delta_{B}(\lambda)=\det(\lambda I - B)= \lambda^{n} + d_{1}
\lambda^{n-1} + \ldots + d_{n-1} \lambda +d_{n},
$$
where $d_{k} = (-1)^{k}\mbox{\large$\delta$}_{k}(B)$, $1 \leq k \leq n$; see \cite[p.$\,$78]{Gantmaher88}.

\vspace{5pt}
Denote by $\bar{B}(\lambda)$ the adjoint matrix of $\lambda I-B$ and let $B_{0}, B_{1}, \ldots, B_{n-2}, B_{n-1}$ be $n$-square
matrices over $K$ determined by
$$
\; \bar{B}(\lambda)
=\mbox{\rm{adj}}(\lambda I\! -\!B)= \lambda^{n-1} B_{0} +
\lambda^{n-2} B_{1} + \ldots + \lambda B_{n-2} + B_{n-1}.
$$

\vspace{5pt}
Recall that $ \; (\lambda I - B) \bar{B}(\lambda) = \Delta_{B}(\lambda)I = \bar{B} (\lambda) (\lambda I\!-\!B). $

\vspace{5pt}
The recurrence $ \; B_{0} \!=\!I; \ B_{k} \!=\!
B_{k-1}  B + d_{k}I \mbox{\ \ \ for \ \ } 1 \leq k < n, $ follows
from equation $\bar{B}(\lambda)  (\lambda I - B)\!=\!
\Delta_B(\lambda) I$, see \cite[p.$\,$91]{Gantmaher88}.

\break

\section{Some results about sums of principal minors}

\bigskip

In this section we give two results about sums of principal minors.

\vspace{5pt}
\begin{theorem}
\label{lemaspona}
For $B \!\in\!K^{n \times n}$ and $\vec{v} = [v_{1}\;\ldots\;v_{n}]^T\!\in\!K^{n\times 1}$ we have$:$
\begin{equation}
\label{Spona} \mbox{\large $\delta$}_{k}^{i} {\big (} B;
\displaystyle\sum_{j=1}^{n}{\!b_{1j} v_{j}}, \ldots,
\displaystyle\sum_{j=1}^{n}{\!b_{nj} v_{j}} {\big )} +
\mbox{\large $\delta$}_{k+1}^{i} {\big (} B; v_{1}, \ldots, v_{n}
{\big )} = \mbox{\large $\delta$}_{k}(B) v_{i}.
\end{equation}
\end{theorem}

\begin{remark}{\rm
The previous result can be described by$:$
\begin{equation}
\mbox{\large $\delta$}_{k}^{i} {\big (} B; B \vec{v} {\big )} +
\mbox{\large $\delta$}_{k+1}^{i} {\big (} B; \vec{v} {\big )} =
\mbox{\large $\delta$}_{k}(B) v_{i} , \ \ 1 \!\leq\! i
\!\leq\! n,
\end{equation}

\noindent
or simply by the following vector equation$:$

\begin{equation}
\label{Sponaaa} \mbox{\large $\vec{\delta}$}_{k} {\big (} B; B
\vec{v} {\big )} + \mbox{\large $\vec{\delta}$}_{k+1} {\big (} B;
\vec{v} {\big )} = \mbox{\large $\delta$}_{k}(B) \vec{v}.
\phantom{ \ \ 1 \!\leq\! i \!\leq\! n,}
\end{equation}}
\end{remark}

\vspace{5pt}
\noindent
{\it Proof of Theorem \ref{lemaspona}.} Let $\vec{e}_{s}\!\in\!K^{n \times 1}$ denote the column whose only
nonzero entry is $1$ in $s^{th}$ position. We also write $B_{\downarrow s}$ for $s^{th}$ column of the matrix $B$
and $[B]_{\hat{s}}$ for a square matrix of order $n\!-\!1$ obtained from $B$ by deleting its $s^{th}$ column and row.
According to the notation used in {\rm (\ref{lepotani})}, let $[B^i(B_{\downarrow s})]_{\hat{s}}$ stand for the matrix
of order $n\!-\!1$ obtained from $B$ by substituting $s^{th}$ column $B_{\downarrow s}$ in place of $i^{th}$ column,
and then by deleting $s^{th}$ column and $s^{th}$ row of the new matrix. By applying linearity of $\vec{\delta}_{k} {\big (} B; \vec{v} {\big )}$
with respect to $\vec{v}$, we have$:$

\vspace{-10pt}
$$
\begin{array}{rcl}
\mbox{\large $\delta$}_{k}^{i}   {\big (} B; B \vec{v} {\big )} +
\mbox{\large $\delta$}_{k+1}^{i} {\big (} B;   \vec{v} {\big )}
\!&\!\!=\!\!& \!
\mbox{\large $\delta$}_{k}^{i} {\big (} B;
\displaystyle\sum\limits_{s=1}^{n} {\! v_s B_{\downarrow s}} {\big
)} \!+\! \mbox{\large $\delta$}_{k+1}^{i} {\big (} B;
\displaystyle\sum\limits_{s=1}^{n} {\! v_s \vec{e}_{s}}
{\big )}=                                                                                           \\[2.0 ex]
\displaystyle\sum\limits_{s=1}^{n} {\! v_s \mbox{\large
$\delta$}_{k}^{i} {\big (} B; B_{\downarrow s} {\big )}} \!+\!
\displaystyle\sum\limits_{s=1}^{n} {\!v_s \mbox{\large
$\delta$}_{k+1}^{i} {\big (} B;  \vec{e}_{s} {\big )}}
\!&\!\!=\!\!&\!
\displaystyle\sum\limits_{s=1}^{n} {\! v_s {\big (} \mbox{\large
$\delta$}_{k}^{i}   {\big (} B; B_{\downarrow s} {\big )} \!+\!
\mbox{\large $\delta$}_{k+1}^{i} {\big (} B;  \vec{e}_{s}    {\big )}
{\big )}}=                                                                                          \\[2.0 ex]
 \;\; v_i {\big (}\, \mbox{\large $\delta$}_{k}^{i}   {\big
(} B; B_{\downarrow i} {\big )} \;+\; \mbox{\large
$\delta$}_{k+1}^{i} {\big (} B;  \vec{e}_{i}     {\big )} \,{\big )}
\!&\!\!+\!\!& \!
\hspace{-2pt} \displaystyle\sum\limits_
{\begin{array}{c}
                                         \\ [-3.75 ex]
\mbox{\scriptsize $s   \!\!=\!\!  1$}    \\ [-1.0 ex]
\mbox{\scriptsize $s  \!\!\neq\!\! i$}   \\
\end{array}}^{n}
{\hspace{-4pt}  v_s {\big (} \mbox{\large $\delta$}_{k}^{i} {\big
(} B; B_{\downarrow s} {\big )} + \mbox{\large $\delta$}_{k+1}^{i}
{\big (} B;  \vec{e}_{s}     {\big )} {\big )}}.
\end{array}
$$

\vspace{-5pt}
\noindent
First, we compute $ v_i {\big (} \mbox{\large $\delta$}_{k}^{i}  {\big (} B; B_{\downarrow i}
{\big )} \!+\! \mbox{\large $\delta$}_{k+1}^{i} {\big (} B; \vec{e}_{i} {\big )} {\big )}
\!=\! v_i {\big (} \mbox{\large $\delta$}_{k}^{i} {\big (} B {\big )} \!+\! \mbox{\large $\delta$}_{k}
{\big (} [B]_{\hat{i}} {\big )}{\big )} \!=\! v_i \mbox{\large $\delta$}_{k} {\big (} B {\big )}.$

\noindent
Then, it remains to show that $\displaystyle\sum\limits_ {
\begin{array}{c}
                                         \\ [-3.75 ex]
\mbox{\scriptsize $s   \!\!=\!\!  1$}    \\ [-1.0 ex]
\mbox{\scriptsize $s  \!\!\neq\!\! i$}   \\
\end{array}}^{n}
{\hspace{-4pt} v_s {\big (} \mbox{\large $\delta$}_{k}^{i}   {\big
(} B; B_{\downarrow s} {\big )} \!+\! \mbox{\large
$\delta$}_{k+1}^{i} {\big (} B;  \vec{e}_{s}     {\big )} {\big
)}}\!=\!0 $.

\vspace{-3pt}
\noindent
It suffices to prove that each term in the sum is zero, i.e.
$$
\mbox{\large $\delta$}_{k}^{i}   {\big (} B;
B_{\downarrow s} {\big )} \!+\! \mbox{\large $\delta$}_{k+1}^{i}
{\big (} B;  \vec{e}_{s}     {\big )} \!=\! 0, \; \mbox{for} \;
s\neq i.
$$
Suppose $s\neq i$. We now consider minors in the sum $\mbox{\large $\delta$}_{k}^{i} {\big (} B; B_{\downarrow s} {\big )}$.
All of them containing $s^{th}$ column are equal to zero, so we deduce
$$
\mbox{\large $\delta$}_{k}^{i} {\big (} B; B_{\downarrow s} {\big
)} \!=\! \mbox{\large $\delta$}_{k}^{i} {\big (} B^i
(B_{\downarrow s}) {\big )} \!=\! \mbox{\large $\delta$}_{k}^{i}
{\big (} [B^i (B_{\downarrow s})]_{\hat{s}} {\big )}.
$$
If $s\neq i$, then each minor in the sum $\mbox{\large $\delta$}_{k+1}^{i} {\big (} B;  \vec{e}_{s} {\big )}$ necessarily
contains $s^{th}$ row and $i^{th}$ column. By interchanging $i^{th}$ and $s^{th}$ column we multiply each minor by $-1$.
We now proceed by expanding these minors along $s^{th}$ column to get $-1$ times the corresponding $k^{th}$ order principal minors
of matrix $B^i (B_{\downarrow s})$ which do not include $s^{th}$ column.
Hence $ \mbox{\large $\delta$}_{k+1}^{i} {\big (} B;  \vec{e}_{s} {\big )}
\!=\! -\mbox{\large $\delta$}_{k}^{i}   {\big (} [B^i
(B_{\downarrow s})]_{\hat{s}} {\big )} $ and the proof is
complete. $\; \Box$

\vspace{5pt}
In the following theorem we give some interesting correspondence between the coefficients $B_{k}$ of the matrix
polynomial $\bar{B}(\lambda)=\mbox{\rm{adj}}(\lambda I\! -\!B)$ and sums of principal minors $\mbox{\large
$\vec{\delta}$}_{k+1}{\big (} B; \vec{v} {\big )}$, $0 \!\leq\! k \!<\! n$.

\vspace{5pt}
\begin{theorem}
\label{lemalepi} Given an arbitrary column
$[v_{1}\;\ldots\;v_{n}]^T\!\in\!K^{n\times 1}$, it holds$:$
\begin{equation}
\label{lepi} B_k \vec{v} \!=\! B_{k} \left[\!
\begin{array}{c}
v_{1}  \\[1.0 ex]
v_{2}  \\[1.0 ex]
\vdots \\[1.0 ex]
v_{n}
\end{array}
\!\right] = (-1)^{k} \!\left[\!
\begin{array}{l}
\mbox{\large $\delta$}_{k+1}^{1}(B; v_{1},\ldots,v_{n})  \\[1.0 ex]
\mbox{\large $\delta$}_{k+1}^{2}(B; v_{1},\ldots,v_{n})  \\[1.0 ex]
\vdots                                                   \\[1.0 ex]
\mbox{\large $\delta$}_{k+1}^{n}(B; v_{1},\ldots,v_{n})
\end{array}
\!\right] =(-1)^{k}\mbox{\large $\vec{\delta}$}_{k+1}{\big (} B;
\vec{v} {\big )}.
\end{equation}
\end{theorem}

\noindent
{\it Proof.} The proof proceeds by induction on $k$. It is being obvious for $k=0$. \\
Assume, as induction hypothesis (IH), that the statement is true
for $k\!-\!1$. Multiplying the right side of the~equation
$B_{k}=B_{k-1} B + d_{k}I$, by the vector $\vec{v}$, we obtain:
$$
\begin{array}{rcl}
B_{k} \vec{v} &\!=\!& B_{k-1} (B \vec{v}) \!+\! d_{k} \vec{v}
\mathop{=}\limits_{(\mbox{\footnotesize IH})}
(-1)^{k-1}\mbox{\large $\vec{\delta}$}_{k}
{\big (} B; B \vec{v} {\big )} \!+\! d_{k} \vec{v}         \\[1.5 ex]
&\!=\!& (-1)^{k-1}{\big (}\mbox{\large $\vec{\delta}$}_{k} {\big
(} B; B \vec{v} {\big )} \!-\! \mbox{\large $\delta$}_{k} \vec{v}
{\big )} \mathop{=}\limits_{(\ref{Sponaaa})} (-1)^{k}\mbox{\large
$\vec{\delta}$}_{k+1} {\big (} B; \vec{v} {\big )}.\;\Box
\end{array}
$$
\begin{remark} {\rm
Theorem {\rm \ref{lemalepi}} seems to have an independent
application. Taking $\vec{v}\!=\!\vec{e}_{j}, 1\leq j\leq n$, we
prove formulae $(8)-(10)$ given in {\rm \cite{Downs75}}.}
\end{remark}

\break

\section{Formulae of total reduction}

\bigskip

We can now obtain a new type of reduction for the linear systems with constant coefficients from \cite{MalesevicTodoricJovovicTelebakovic2010}
applying results of previous section. For the sake of completeness we introduce some definition.

Let $K$ be a field, $V$ a vector space over field $K$ and let $A:V\longrightarrow V$ be a linear operator. We consider a linear
system of first order $A$-operator equations with constant coefficients in unknowns $x_{i}$, $1 \!\leq\! i \!\leq\! n$,~is$:$
\begin{equation}
\label{B-Sistem-1} \left\{
\begin{array}{rcl}
A(x_1) & \!\!=\!\!      & b_{11} x_1 + b_{12} x_2 + \ldots + b_{1n} x_n + \varphi_{1} \\[0.5 ex]
A(x_2) & \!\!=\!\!      & b_{21} x_1 + b_{22} x_2 + \ldots + b_{2n} x_n + \varphi_{2} \\[0.5 ex]
       & \!\!\vdots\!\! &                                                             \\[0.5 ex]
A(x_n) & \!\!=\!\!      & b_{n1} x_1 + b_{n2} x_2 + \ldots +
b_{nn} x_n + \varphi_{n}
\end{array}
\right\}
\end{equation}
for $b_{ij}\!\in\!K$ and $\varphi_{i}\!\in\!V$. We say that $B\!=\![b_{ij}]_{i, j=1}^n\!\in\!K^{n \times n}$ is the system matrix,
and $\vec{\varphi}\!=\![\varphi_{1}\;\ldots\;\varphi_{n}]^T\!\in\!V^{n\times 1}$ is the free column.

\vspace{5pt}
Let $\vec{x}\!=\![x_{1}\;\ldots\;x_{n}]^T$ be a column of unknowns
and $\vec{A} : V^{n\times 1}\longrightarrow V^{n\times 1}$ be a
vector operator defined componentwise
$\vec{A}(\vec{x})\!=\![A(x_1)\;\ldots\;A(x_n)]^T.$ Then system
(\ref{B-Sistem-1}) can be written in the following vector form$:$
\begin{equation}
\label{B-Sistem-2} \hspace{10pt}
\vec{A}(\vec{x})=B \vec{x}+\vec{\varphi}.
\end{equation}
Any column $\vec{v}\!\in\! V^{n \times 1}$ which satisfies the previous system is its solution.

\vspace{5pt}
Powers of operator $A$ are defined as usual $A^{i}\!=\!A^{i-1}\!\circ A$ assuming that $A^{0}\!=\!I~:~V~\longrightarrow~V$
is the identity operator. By $n^{th}$ order linear $A$-operator equation with constant coefficients, in unknown $x$,
we mean$:$
\begin{equation}
\label{Eq_5}
A^{n}(x) + b_{1} A^{n-1}(x) + \ldots + b_{n} I(x) =\varphi,
\end{equation}
where $b_1,\! \ldots\!, b_n \!\in\! K$ are coefficients and $\varphi \in V$. Any vector $v\in V$ which satisfies ({\ref{Eq_5}) is it's solution.

\vspace{5pt}
The following theorem separates variables of the initial system.

\begin{theorem}
\label{Teorema TRS 1} Assume that the linear system of first order $A$-operator equations is given in the form
{\rm (\ref{B-Sistem-2})}, $\vec{A}(\vec{x}) = B \vec{x}+\vec{\varphi}$, and that matrices $B_{0}, \ldots, B_{n-1}$
are coefficients of the matrix polynomial $\bar{B}(\lambda)=\mbox{\rm{adj}}(\lambda I\!-\!B)$. Then it holds$:$
\begin{equation}
\label{TRS 1}
{\big (}\Delta_B{\big (}\vec{A}{\big )}{\big )}(\vec{x})
\!= \!
\displaystyle\sum\limits_{k=1}^{n}{B_{k-1}\vec{A}^{n-k}(\vec{\varphi})}.
\end{equation}
\end{theorem}

\noindent
{\it Proof.}
Let $L_B : V^{n \times 1}\longrightarrow V^{n \times 1}$ be a linear operator defined by $L_B(\vec{x})=B \vec{x}$.
Replacing $\lambda I$ by $\vec{A}$ in the equation $\Delta_B(\lambda) I = \bar{B}(\lambda) (\lambda I - B)$ we obtain

\bigskip

\noindent
$$
\begin{array}{rclllcl}
  \Delta_B{\big (}\vec{A}\,{\big )}
  & \!\!=\!\! &
  \bar{B}{\big (}\vec{A}\,{\big )}\circ {\big (}\vec{A}-L_B{\big )},
                                                                     & \!\!\!\!\mbox{hence}\!\!\!\! &&&      \\[2.5 ex]
  \Delta_B{\big (}\vec{A}\,{\big )}(\vec{x})
  & \!\!=\!\! &
  \bar{B}{\big (}\vec{A}\,{\big )}{\big (}{\big (}\vec{A}-L_B{\big )}(\vec{x}){\big )}             &&&&      \\[0.5 ex]
  & \!\!=\!\! &
  \bar{B}{\big (}\vec{A}\,{\big )}{\big (}\vec{A}(\vec{x})-B \vec{x}{\big )}
  &\!\!\!\!\!\!\mathop{=}\limits_{\mbox{\footnotesize (\ref{B-Sistem-2})}}\;\;\;\;\!\!&
 \!\!\!\!\!\!\!\!\bar{B}{\big (}\vec{A}\,{\big )}(\vec{\varphi})
  &\!\!\!\! = \!\!\!\!&
  \displaystyle\sum\limits_{k=1}^{n}{B_{k-1}\vec{A}^{n-k}(\vec{\varphi})}. \;\; \Box
\end{array}
$$

\vspace{5pt}
The next theorem is an operator generalization of Cramer's rule.

\vspace{5pt}
\begin{theorem}
\label{Teorema TRS 2} $($The theorem of total reduction - vector form$)$ \\
Linear system of first order $A$-operator equations {\rm (\ref{B-Sistem-2})}
can be reduced to the system of $n^{th}$ order $A$-operator equations$:$
\begin{equation}
\label{TRS} {\big (} \Delta_B{\big (} \vec{A}\, {\big )} {\big )}({\vec{x}})
\! = \!
\displaystyle\sum\limits_{k=1}^{n}(-1)^{k-1} \mbox{\large $\vec{\delta}$}_{k} {\big (} B;{\vec{A}}^{n-k}({\vec{\varphi}}){\big )}.
\end{equation}
\end{theorem}

\noindent
{\it Proof.} It is an immediate consequence of Theorem \ref{Teorema TRS 1} and Theorem \ref{lemalepi}$:$
$$
\begin{array}{rclcl}
\Delta_B{\big (} \vec{A}\, {\big )}({\vec{x}}) &\!
\mathop{=}\limits_{(\ref{TRS 1})} \!&
\displaystyle\sum\limits_{k=1}^{n}
B_{k-1}\vec{A}^{n-k}(\vec{\varphi}) &\!
\mathop{=}\limits_{(\ref{lepi})} \!&
\displaystyle\sum\limits_{k=1}^{n} (-1)^{k-1} \mbox{\large
$\vec{\delta}$}_{k} {\big (} B; \vec{A}^{n-k}(\vec{\varphi}){\big
)}. \;\Box
\end{array}
$$

\vspace{5pt}
We can now rephrase the previous theorem as follows.

\vspace{5pt}
\begin{theorem}
\label{Teorema TRS 3} $($The theorem of total reduction$)$  \\
Linear system of first order $A$-operator equations
{\rm(\ref{B-Sistem-1})} implies the system, which consists of $n^{th}$ order $A$-operator equations$:$
\begin{equation}
\label{TRS 3} \left\{
\begin{array}{rcl}
{\big (}\Delta_B (A){\big )}(x_1) &\! = \!&
\displaystyle\sum\limits_{k=1}^{n} (-1)^{k-1} \mbox{\large
$\delta$}_{k}^{1} {\big (} B; \vec{A}^{n-k}(\vec{\varphi}){\big )}
\\ [-1.0 ex] & \vdots &
\\ [-1.0 ex] {\big (}\Delta_B (A){\big )}(x_i) &\! = \!&
\displaystyle\sum\limits_{k=1}^{n} (-1)^{k-1} \mbox{\large
$\delta$}_{k}^{i} {\big (} B; \vec{A}^{n-k}(\vec{\varphi}){\big )}
\\ [-1.0 ex] & \vdots &
\\ [-1.0 ex] {\big (}\Delta_B (A){\big )}(x_n) &\! = \!&
\displaystyle\sum\limits_{k=1}^{n} (-1)^{k-1} \mbox{\large
$\delta$}_{k}^n {\big (} B; \vec{A}^{n-k}(\vec{\varphi}){\big )}
\end{array}
\right\}.
\end{equation}
\end{theorem}

\vspace{5pt}
\begin{remark}
{\rm System {\rm (\ref{TRS 3})} has separated variables and it is called totally reduced. The obtained system is suitable for applications since
it does not require a change~of~base. This system consists of $n^{th}$ order linear $A$-operator equations~which differ only in the variables and
in the non-homogeneous terms.}
\end{remark}

\vspace*{2.5 pt}
Transformations of the linear systems of operator equations into independent equations are important in applied mathematics \cite{Fortmann--Hitz--77}.
In the following two sections we apply our theorem of total reduction to the specific linear operators~$A$.

\medskip

\section{Cauchy problem}

\bigskip

Let us {assume that} $A$ {is a} differential operator on the vector space of real functions and {that} system~{\rm(\ref{B-Sistem-1})} has initial
conditions $x_{i}(t_{0}) \!=\! c_{i},$ for \mbox{$1\!\leq\!i\!\leq\!n$}. {Then} the Cauchy problem for system {\rm(\ref{B-Sistem-1})} has a unique solution.
\mbox{Using} form (\ref{B-Sistem-2}), we obtain additional $n\!-\!1$ initial conditions of  $i^{th}$ equation in system~{\rm(\ref{TRS 3})}$:$
\begin{equation}
\label{DJ_Pocetni_Uslovi}
\begin{array}{rcl}
{\big (}A^{j}(x_{i}){\big )}(t_{0}) &\!\!\! = \!\!\!\!& {\big
[}B^{j} \vec{x}(t_{0}){\big ]}_{i} \!+\!
\displaystyle\sum\limits_{k=0}^{j-1}{\big [}B^{j-1-k}
\vec{A}^{k}(\vec{\varphi})(t_{0}){\big ]}_{i}, \;\; 1 \leq j \leq n\!-\!1,
\end{array}
\end{equation}
where $[\,\,]_{i}$ denotes $i^{th}$ coordinate. Then each equation
in system {\rm(\ref{TRS 3})} has a unique solution under given
conditions and additional conditions (\ref{DJ_Pocetni_Uslovi}),
and these solutions form a unique solution to system
{\rm(\ref{B-Sistem-1})}. Therefore, formulae {\rm(\ref{TRS 3})}
can be used for solving systems of differential equations.

\vspace{2.5 pt}
It is worth pointing out that the above method can be also extended to systems of difference equations.

\medskip

\section{Differential transcendency}

\bigskip

Now suppose that $V$ is {the} vector space of meromorphic functions over the complex field $C$ and that $A$ is
a differential operator, $A(x)\!=\!\mbox{\large $\frac{d}{d\,z}$}(x)$.
Let us consider system {\rm (\ref{B-Sistem-1})} under these assumptions.

\vspace{2.5 pt}
Recall that a function $x \in V$ is differentially algebraic if it satisfies a differential algebraic equation with
coefficients {from} $C$, otherwise {it} is differentially transcendental (see
\cite{Markus07}--\cite{Cruz-Victoria--Martínez-Guerra--Rincón-Pasaye--08}).

\vspace{2.5 pt}
Let us consider non-homogenous linear differential equation of $n^{th}$ order in the form {\rm (\ref{Eq_5})} where
$b_{1}, \ldots, b_{n} \in C$ are constants and $\varphi \in V$. If $x$ is differentially transcendental then
$\Delta_B (A)(x)$ is {also} a differentially transcendental function. On the other hand, if $\varphi$ is differentially transcendental,
then, based~on Theorem 2.8. from \cite{MijajlovicMalesevic08}, the solution $x$ of equation {\rm (\ref{Eq_5})}
is a differentially transcendental function. Therefore we obtain the equivalence.
\begin{theorem}
\label{Th_DT}
Let $x$ be a solution of the equation {\rm (\ref{Eq_5})}, then $x$ is a differentially transcendental function iff $\varphi$ is
a differentially transcendental function.
\end{theorem}

We also have the following statement about differential transcendency.
\begin{theorem}
Let $\varphi_{j}$ be the only differentially transcendental component of the free column $\vec{\varphi}$. Then for any solution $\vec{x}$
of system {\rm (\ref{B-Sistem-2})}, the corresponding entry $x_{j}$ is also {a} differentially transcendental function.
\end{theorem}
{\it Proof.}
The sum $\sum_{k=1}^{n}{(-1)^{k-1} \mbox{\large $\delta$}_{k}^{j}{\big (} B; \vec{A}^{n-k}(\vec{\varphi}){\big )}}$
must be a differentially transcendental function. {The} previous theorem {applied}
to the equation$:$
\begin{equation}
\hspace{-45pt} \Delta_B (A)(x_j) = \sum_{k=1}^{n}{(-1)^{k-1}
\mbox{\large $\delta$}_{k}^{j}{\big (} B;
\vec{A}^{n-k}(\vec{\varphi}){\big )}},
\end{equation}
implies that $x_{j}$ is a differentially transcendental function too. $\;\Box$

\vspace{7.5 pt}
Let us consider system {\rm (\ref{B-Sistem-1})} and let $\varphi_{1}$ be the only differentially transcendental component of the free column $\vec{\varphi}$.
Then, the coordinate $x_{1}$ is a differentially transcendental function too.
Whether the other coordinates $x_{k}$ are differentially algebraic depends on the system matrix $B$.
From the Formulae of total reduction and Theorem~\ref{Th_DT}. we obtain the following statement.
\begin{theorem}
Let $\varphi_1$ be the only differentially transcendental component of the free column $\vec{\varphi}$ of system {\rm (\ref{B-Sistem-1})}.
Then the coordinate $x_{k}$, \mbox{$k \neq 1$}, of the solution $\vec{x}$, is differentially algebraic iff
in the sum
$
\sum_{j=1}^{n}{\!(-1)^{j-1} \mbox{\large $\delta$}_{j}^{k}{\big (} B; \vec{A}^{n-j}(\vec{\varphi}){\big )}}
$
{appears} no function~\mbox{$A^{n-j}(\varphi_{1}) \; (j \!=\! 1,\ldots,n)$}.
\end{theorem}
\begin{example}
{\rm Let us consider system {\rm (\ref{B-Sistem-1})} in the form {\rm (\ref{TRS 3})} in dimensions \mbox{$n \!=\! 2$}
and \mbox{$n \!=\! 3$} with $\varphi_1$ as the only differentially transcendental component.
The function $x_{1}$ is differentially transcendental.
For~$n \!=\! 2$ the function $x_2$ is differentially algebraic iff \mbox{$b_{21} \!=\! 0$}. For~$n \!=\! 3$ the function $x_2$ is differentially algebraic iff
\mbox{$b_{21} \!=\! 0 \wedge b_{31} \!\cdot\! b_{23} \!=\! 0$} and the function $x_3$ is differentially algebraic iff
\mbox{$b_{31} \!=\! 0 \wedge b_{31} \!\cdot\! b_{22} \!=\! 0$}. $\;\Box$}
\end{example}
Let us emphasize {that} if we consider two or more differentially transcendental components of the free column $\vec{\varphi}$, then {the} differential transcendency of the solution coordinates also depends on some kind of their differential independence
(see for example \cite{Markus07}). 

\bigskip

\noindent
\mbox{\textsc{Acknowledgment.}} Research is partially supported by the Ministry of Science and Education of the Republic of Serbia,
Grant No. ON 174032.

\break

\end{document}